\documentclass{amsart}

\usepackage{amssymb}
\usepackage{mathrsfs}
\usepackage{amscd}
\usepackage{verbatim}

\usepackage[colorlinks,linkcolor={blue},
citecolor={blue},urlcolor={red},]{hyperref}

\theoremstyle{plain}
\newtheorem{theorem}{Theorem}[section]
\theoremstyle{remark}

\newtheorem{example}[theorem]{Example}
\theoremstyle{plain}
\newtheorem{corollary}[theorem]{Corollary}
\newtheorem{lemma}[theorem]{Lemma}
\newtheorem{proposition}[theorem]{Proposition}

\numberwithin{equation}{section}

\def\Z{{\mathbb Z}}
\def\Q{{\mathbb Q}}
\def\R{{\mathbb R}}

\newcommand{\E}{{\mathbb E}}
\renewcommand{\P}{{\mathbb P}}
\newcommand{\F}{{\mathcal F}}
\newcommand{\G}{{\mathcal G}}
\renewcommand{\H}{{\mathcal H}}

\renewcommand{\a}{\alpha}

\newcommand{\g}{\gamma}
\renewcommand{\d}{\delta}
\newcommand{\e}{\varepsilon}

\renewcommand{\O}{\Omega}

\renewcommand{\L}{L^2(0,T)}
\newcommand{\LH}{L^2(0,T;H)}
\renewcommand{\gg}{\g(\L,E)}
\newcommand{\ggH}{\g(\LH,E)}

\newcommand{\beq}{\begin{equation}}
\newcommand{\eeq}{\end{equation}}
\newcommand{\bal}{\begin{aligned}}
\newcommand{\eal}{\end{aligned}}
\newcommand{\ben}{\begin{enumerate}}
\newcommand{\een}{\end{enumerate}}
\newcommand{\bit}{\begin{itemize}}
\newcommand{\eit}{\end{itemize}}

\newcommand{\bth}{\begin{theorem}}
\renewcommand{\eth}{\end{theorem}}
\newcommand{\bpr}{\begin{proposition}}
\newcommand{\epr}{\end{proposition}}
\newcommand{\ble}{\begin{lemma}}
\newcommand{\ele}{\end{lemma}}
\newcommand{\bpf}{\begin{proof}}
\newcommand{\epf}{\end{proof}}
\newcommand{\bex}{\begin{example}}
\newcommand{\eex}{\end{example}}
\newcommand{\bre}{\begin{example}}
\newcommand{\ere}{\end{example}}

\newcommand{\D}{{\mathcal D}}
\newcommand{\calL}{{\mathcal L}}
\newcommand{\n}{\Vert}
\newcommand{\one}{{{\bf 1}}}
\newcommand{\embed}{\hookrightarrow}
\newcommand{\s}{^*}
\newcommand{\lb}{\langle}
\newcommand{\rb}{\rangle}

\newcommand{\limn}{\lim_{n\to\infty}}
\newcommand{\limk}{\lim_{k\to\infty}}
\newcommand{\limj}{\lim_{j\to\infty}}

\newcommand{\bbF}{{\mathbb{F}}}

\newcommand{\wh}{\widehat}

\begin{document}

\title[Stochastic integrability]
{Conditions for stochastic integrability in UMD Banach spaces}

\author{Jan van Neerven}
\address{Delft Institute of Applied Mathematics\\
Delft University of Technology \\ P.O. Box 5031\\ 2600 GA Delft\\The
Netherlands} \email{J.M.A.M.vanNeerven@tudelft.nl}

\author{Mark Veraar}
\address{Delft Institute of Applied Mathematics\\
Delft University of Technology \\ P.O. Box 5031\\ 2600 GA Delft\\The
Netherlands} \email{M.C.Veraar@tudelft.nl}

\author{Lutz Weis}
\address{Mathematisches\, Institut\, I \\
\, Technische\, Universit\"at \, Karlsruhe \\
\, D-76128 \, Karls\-ruhe\\Germany}
\email{Lutz.Weis@math.uni-karlsruhe.de}

\thanks{The first and second named authors are supported by
a `VIDI subsidie' (639.032.201) in the `Vernieuwingsimpuls'
programme of the Netherlands Organization for Scientific Research
(NWO). The first named author is also supported by the Research
Training Network ``Evolution Equations for Deterministic and
Stochastic Systems'' (HPRN-CT-2002-00281). The second named author
gratefully acknowledges support from the Marie Curie Fellowship
Program for a stay at the TU Karlsruhe. The third named author is
supported by grants from the Volkswagenstiftung (I/78593) and the
Deutsche Forschungsgemeinschaft (We 2847/1-2).}



\date{\today}

\begin{abstract}
A detailed theory of stochastic integration in UMD Banach spaces has
been developed recently in \cite{NVW}. The present paper is aimed at
giving various sufficient conditions for stochastic integrability.
\end{abstract}

\maketitle

\section{Introduction}
In the paper \cite{NVW} we developed a detailed theory of stochastic
integration in UMD Banach spaces and a number of necessary and
sufficient conditions for stochastic integrability of processes with
values in a UMD space were obtained. The purpose of the present
paper is to complement these results by giving further conditions
for stochastic integrability.

In Section \ref{sec:approx}, we prove a result announced in
\cite{NVW} on the strong approximation of stochastically integrable
processes by elementary adapted processes. In Section
\ref{sec:domination} we prove two domination results. In Section
\ref{sec:besov} we state a criterium for stochastic integrability in
terms of the smoothness of the trajectories of the process. This
criterium is based on a recent embedding result due to Kalton and
the authors \cite{KNVW}. In Section \ref{sec:besov2} we give an
alternative proof of a special case of the embedding result from
\cite{KNVW} and we prove a converse result which was left open
there. In the final Section \ref{sec:BFS} we give square function
conditions for stochastic integrability of processes with values in
a Banach function spaces.

We follow the notations and terminology of the paper \cite{NVW}.

\section{Approximation}\label{sec:approx}
Throughout this note, $(\O,\F,\P)$ is a probability space endowed
with a filtration $\bbF=(\F_t)_{t\in [0,T]}$ satisfying the usual
conditions, $H$ is a separable real Hilbert space with inner product
$[\cdot,\cdot]_H$, and $E$ is a real Banach space with norm
$\n\cdot\n$. The dual of $E$ is denoted by $E\s$.

We call an operator-valued stochastic process
$\Phi:[0,T]\times\Omega\to \calL(H,E)$ {\em elementary adapted} with
respect to the filtration $\bbF$ if it is of the form
\[
\Phi = \sum_{n=1}^N \sum_{m=1}^M \one_{(t_{n-1},t_n]\times A_{mn}}
\sum_{k=1}^K h_k \otimes x_{kmn},
\]
where $0\le t_0 < \dots< t_N \le T$ with the convention that
$(t_{-1},t_0] = \{0\}$,  the sets $A_{1n},\dots,A_{Mn}\in
\F_{t_{n-1}}$ are
 disjoint for all $n=1,\dots,N$, and the vectors
$h_1,\dots,h_K \in H$ are orthonormal.
% [[[ if the index n starts at 0:
%Here we use the conventions $(t_{-1},t_0] = \{0\}$ and $\F_{t_{-1}} = \F_0$.

Let $W_H = (W_H(t))_{t\in [0,T]}$ be an {\em $H$-cylindrical
Brownian motion}, i.e., each $W_H(t)$ is a bounded operator from $H$
to $L^2(\O)$, for all $h\in H$ the process $W_H h = (W_H(t)h)_{t\in
[0,T]}$ is a Brownian motion, and for all $t_1,t_2\in [0,T]$ and
$h_1,h_2\in H$ we have
\[\E(W_H(t_1)h_1 \, W_H(t_2)h_2)
= (t_1\wedge t_2)[h_1,h_2]_H.\] We will always assume that $W_H$ is
adapted to $\bbF$, i.e., each Brownian motion $W_H h$ is adapted to
$\bbF$. The stochastic integral of an elementary adapted process
$\Phi$ of the above form with respect to $W_H$ is defined in the
obvious way as
\[\int_0^t \Phi\,dW_H = \sum_{n=1}^N \sum_{m=1}^M
\one_{A_{mn}}\sum_{k=1}^K \big(W_H(t_n\wedge t)h_k -
W_H(t_{n-1}\wedge t)h_k\big) \otimes x_{kmn}.
\]

A process $\Phi:[0,T]\times\O\to\calL(H,E)$ is called {\em
$H$-strongly measurable} if for all $h\in H$, $\Phi h$ is strongly
measurable. Similarly, $\Phi$ is {\em $H$-strongly adapted} if for
all $h\in H$, $\Phi h$ is strongly adapted.

An $H$ strongly measurable and adapted process
$\Phi:[0,T]\times\O\to\calL(H,E)$ is called {\em stochastically
integrable} with respect to $W_H$ if there exists a sequence of
elementary adapted processes $\Phi_n:[0,T]\times\O\to\calL(H,E)$ and
a $\zeta:\O\to C([0,T];E)$ such that \ben
\item[(i)]
$\displaystyle\limn\lb\Phi_n h,x\s\rb=\lb\Phi h,x\s\rb$ in measure
%on $[0,T]\times\O$
for all $h\in H$ and $x\s\in E\s$;
\item[(ii)]
$\displaystyle \limn \int_0^\cdot\!\Phi_n\,dW_H = \zeta$ measure in
$C([0,T];E)$. \een The process $\zeta$ is uniquely determined almost
surely. We call $\zeta$ the {\em stochastic integral} of $\Phi$,
notation:
\[ \zeta =: \int_0^\cdot \Phi\,dW_H.\]

It is an easy consequence of (i), (ii), and \cite[Proposition
17.6]{Kal} that if $\Phi$ is stochastically integrable, then for all
$x\s\in E\s$ we have
\[\limn \Phi_n\s x\s =\Phi\s x\s \ \ \hbox{in $\LH$ almost surely.}\]
%Indeed, by BDG the processes
%$\Phi_n\s x\s: [0,T]\times\O\to H$ define a Cauchy sequence in $L^p(\O;\LH)$. %Let $\phi_{x\s}$ be its limit.  Then by (i), for all $h\in H$ and $x\s\in %E\s$ we have
%$[\phi_{x\s}, h]_H = [\Phi\s x\s, h]_H$ in $L^p(\O;\L)$.
%It follows that $\Phi\s x\s\in L^p(\O;\LH)$ and $$\limn \Phi_n\s x\s =\Phi\s %x\s \ \ \hbox{in $L^p(\O;\LH)$.}$$

For UMD spaces $E$ we show that in the definition of stochastic
integrability it is possible to strengthen the convergence of the
processes $\Phi_n h$ in (i) to strong convergence in measure. The
main result of this section was announced without proof in
\cite{NVW} and is closely related to a question raised by McConnell
\cite[page 290]{MCC}.

\begin{theorem}\label{thm:main} Let $E$ be a UMD space.
If the process $\Phi:[0,T]\times\O\to \calL(H,E)$ is $H$-strongly
measurable and adapted and stochastically integrable with respect to
$W_H$, there exists a sequence of elementary adapted processes
$\Phi_n:[0,T]\times\O\to \calL(H,E)$ such that \ben
\item[(i)$'$]
$\displaystyle\limn  \Phi_n h= \Phi h$ in measure
%on $[0,T]\times\O$
for all $h\in H$;
\item[(ii)]
$\displaystyle\limn \int_0^\cdot \Phi_n\,dW_H =\int_0^\cdot
\Phi\,dW_H $ in measure in $C([0,T];E)$. \een
\end{theorem}

For the definition of the class of UMD Banach spaces and some of its
applications in Analysis we refer to Burkholder's review article
\cite{Bu}.

\medskip
Let $\H$ be a separable real Hilbert space and let $(g_n)_{n\ge 1}$
be a sequence of independent standard Gaussian random variables on a
probability space $(\O',\F',\P')$. A linear operator $R: \H\to E$ is
said to be {\em $\g$-radonifying} if for some (every) orthonormal
basis $(h_n)_{n\ge 1}$ of $\H$ the Gaussian sum
 $\sum_{n\ge 1} g_n \, Rh_n$ converges in $L^2(\O';E)$.
The linear space of all $\g$-radonifying operators from $\H$ to $E$
is denoted by $\g(\H,E)$. This is space is a Banach space endowed
with the norm
\[ \n R
\n_{\g(\H,E)} := \Big(\E' \Big\n\sum_{n\ge 1} g_n \,
Rh_n\Big\n^2\Bigr)^\frac12 .\] For more information we refer to
\cite{Bo,DJT,KW}. The importance of spaces of $\g$-radonifying
operators in the theory of stochastic integration in infinite
dimensions is well established; see \cite{NVW,NW} and the references
given therein.

An $H$-strongly measurable function $\Phi:[0,T]\to \calL(H,E)$ is
said to {\em represent} an element $R\in \ggH$ if for all $x\s\in
E\s$ we have $\Phi^* x^*\in \LH$ and, for all $f\in\LH$,
\[
\lb Rf, x\s\rb = \int_0^T [f(t), \Phi^*(t) x\s]\,dt.
\]

Extending the above definition, we say that an $H$-strongly
measurable process $\Phi:[0,T]\times \O\to \calL(H,E)$ {\em
represents} a random variable $X: \O\to \ggH$ if for all $x\s\in
E\s$ almost surely we have $\Phi\s x\s\in \LH$ and, for all
$f\in\LH$,
\[
\lb X f, x\s\rb = \int_0^T [f(t), \Phi\s(t)x\s]_H\,dt \ \
\text{almost surely}.
\]

For $H$-strongly measurable process we have the following simple
result \cite[Lemma 2.7]{NVW}.

\begin{lemma}\label{lem:ae-repr}
Let $\Phi:[0,T]\times\O\to \calL(H,E)$ be an $H$-strongly measurable
process and let  $X:\O\to\ggH$ be strongly measurable. The following
assertions are equivalent:

\begin{enumerate}
\item $\Phi$ represents $X$.
\item $\Phi(\cdot, \omega)$ represents $X(\omega)$ for almost all $\omega\in \O$.
\end{enumerate}
\end{lemma}
%The assumptions of (2) already imply that the induced mapping $\omega\mapsto
%X(\omega)$ from $\O$ to $\ggH$ has a strongly measurable version (see
%\cite[Remark 2.8]{NVW}).

%The special case for functions (which we regard as `deterministic' processes)
%can be found in \cite[Theorem 4.2]{NW}; see \cite{Ma,Ro,RS} for related
%results. An $H$-strongly measurable process $\Phi$ is said to be {\em adapted}
%if $\Phi h$ is adapted for all $h\in H$.

%We call the process $\Phi$ {\em scalarly measurable} if
%the processes $\lb \Phi h,x\s\rb$ are measurable for all $h\in H$ and $x\s\in E\s$.

For a Banach space $F$ we denote by $L^0(\O;F)$ the vector space of
all $F$-valued random variables on $\O$, identifying random
variables if they agree almost surely. Endowed with the topology of
convergence in measure, $L^0(\O;F)$ is a complete metric space. The
following result is obtained in \cite{NVW}.

\begin{proposition}\label{prop:NVW}
Let $E$ be a UMD space. For an $H$-strongly measurable and adapted
process $\Phi:[0,T]\times \O\to \calL(H,E)$ the following assertions
are equivalent: \ben
\item $\Phi$ is stochastically integrable with respect to $W_H$;
\item $\Phi$ represents a random variable $X:\O\to \ggH$.
\een
%Moreover, there exist constants $0 < c\le C<\infty$, depending only on $p$ and
%$E$, such that
%$$c^p \n X\n_{L^p(\O;\ggH)}^p \le
%\E \Big\n\int_0^T \Phi\,dW_H\Big\n^p \le C^p \n X\n_{L^p(\O;\ggH)}^p.$$
In this case $\Phi h$ is stochastically integrable with respect to
$W_H h$ for all $h\in H$, and for every orthonormal basis
$(h_n)_{n\ge 1}$ of $H$ we have
\[ \int_0^\cdot \Phi\,dW_H = \sum_{n\ge 1}  \int_0^\cdot \Phi h_n\,dW_H
h_n,\] with almost sure unconditional convergence of the series
expansion in $C([0,T];E)$. Moreover, for all $p\in (1, \infty)$
\[
\E\sup_{t\in [0,T]}\Bigl\n\int_0^t \Phi(s)\,dW_H(s)\Bigr\n^p
\eqsim_{p,E} \E\n X\n_{\ggH}^p.
\]
Furthermore, the mapping $X\mapsto \int_0^\cdot \Phi\,dW_H$ has a
unique extension to a continuous mapping
\[
L^0(\O;\ggH)\to L^0(\O;C([0,T];E)).
\]
\end{proposition}

As we will show in a moment, the series expansion in Proposition
\ref{prop:NVW} implies that in order to prove Theorem \ref{thm:main}
it suffices to prove the following weaker version of the theorem:

\begin{theorem}\label{thm:mainE}
Let $E$ be a UMD space. If the process $\phi:[0,T]\times\O\to E$ is
strongly measurable and adapted and stochastically integrable with
respect to a Brownian motion $W$, there exists a sequence of
elementary adapted processes $\phi_n:[0,T]\times\O\to E$ such that
\ben
\item[(i)$'$]
$\displaystyle\limn  \phi_n = \phi$ in measure;
\item[(ii)]
$\displaystyle\limn \int_0^\cdot \phi_n\,dW =\int_0^\cdot \phi\,dW $
in measure in $C([0,T];E)$. \een
\end{theorem}

This theorem may actually be viewed as the special case of Theorem
\ref{thm:main} corresponding to $H=\R$, by identifying $\calL(\R,E)$
with $E$ and identifying $\R$-cylindrical Brownian motions with
real-valued Brownian motions.

\medskip
To see that Theorem \ref{thm:main} follows from Theorem
\ref{thm:mainE} we argue as follows. Choose an orthonormal basis
$(h_n)_{n\ge 1}$ of $H$ and define the processes
$\Psi_n:[0,T]\times\O\to\calL(H,E)$ by
\[ \Psi_n h :=
\sum_{j=1}^n [h,h_j]_H\, \Phi h_j.
\]
Clearly, $\limn  \Psi_n h= \Phi h$ pointwise, hence in measure, for
all $h\in H$. In view of the identity
\[ \int_0^\cdot \Psi_n\,dW_H =
\sum_{j=1}^n\int_0^\cdot \Phi h_j \,dW_H h_j\] and the series
expansion in Proposition \ref{prop:NVW}, we also have
\[\limn \int_0^\cdot \Psi_n\,dW_H
=\int_0^\cdot \Phi\,dW_H \ \ \hbox{ in measure in $C([0,T];E)$.}
\]
With Theorem \ref{thm:mainE}, for each $n\ge 1$ we choose a sequence
of elementary adapted processes $\phi_{j,n}: [0,T]\times\O\to E$
such that $\limj \phi_{j,n} = \Phi h_n$ in measure and
\[\limj \int_0^\cdot \phi_{j,n}\,dW_H h_n
=\int_0^\cdot \Phi h_n\,dW_H h_n \ \ \hbox{in measure in
$C([0,T];E)$.}
\]
Given $k\ge 1$, choose $N_k\ge 1$ so large that
%$$ \Big|\Big\{\n \Phi h_n - \Phi_{N_k} h_n\n>\frac1k\Big\} \Big| < %\frac{1}{k}, \qquad n=1,\dots,k,$$ and
\[ \P \Big\{\Big\n \int_0^\cdot \Phi-\Psi_{N_k}\,dW_H\Big\n_\infty>\frac1k\Big\} <
\frac1{k}.
\]
Let $\lambda$ be denoted for the Lebesgue measure on $[0,T]$. For
each $n=1,\dots,N_k$ choose $j_{k,n}\ge 1$ so large that
\[ \lambda\otimes\P\Big\{\n
\Phi h_n - \phi_{j_{k,n},n}\n>\frac1{kN_k}\Big\}  < \frac{1}{kN_k}
\]
and
\[ \P\Big\{\Big\n \int_0^\cdot \Phi h_n - \phi_{j_{k,n},n}\,dW_H
h_n\Big\n_\infty>\frac1k\Big\} < \frac1{kN_k}.
\]
Define $\Phi_k:[0,T]\times\O\to \calL(H,E)$ by
\[
\Phi_k h:=  \sum_{n=1}^{N_k} [h,h_n]_H\, \phi_{j_{k,n},n}, \qquad
h\in H.
\]
Each $\Phi_k$ is elementary adapted. For all $h\in H$ with $\n
h\n_H=1$ and all $\d>0$ we have, for all $k\ge 1/\delta$,
\[\bal \ & |\{\n \Phi h - \Phi_k
h\n>2\d\} |
\\ & \quad\qquad \le
\lambda\otimes\P\{\n \Phi h - \Psi_{N_k} h\n>\d\}  +
\lambda\otimes\P\{\n \Psi_{N_k} h - \Phi_k h\n>\d\}
\\ & \quad\qquad<  \lambda\otimes\P\{\n \Phi h - \Psi_{N_k}h\n>\d\}  + \frac1k.
\eal
 \]
Hence, $\limk \Phi_k h = \Phi h$ in measure for all $h\in H$. Also,
\[
\bal \ & \P\Big\{\Big\n \int_0^T \Phi - \Phi_k\,dW_H\Big\n_\infty
>\frac2k\Big\}
\\ & \quad  \le
 \P\Big\{\Big\n \int_0^T \Phi - \Psi_{N_k}\,dW_H\Big\n_\infty>\frac1k\Big\}
 + \P\Big\{ \Big\n \int_0^T \Psi_{N_k} - \Phi_k\,dW_H\Big\n_\infty>\frac1k\Big\}
\\ & \quad< \frac1k+\frac1k = \frac{2}{k},
\eal
\]
and therefore $\displaystyle \limk\int_0^\cdot
\Phi_k\,dW_H=\int_0^\cdot\Phi\,dW_H$ in measure in $C([0,T];E)$.
Thus the processes $\Phi_k$ have the desired properties.
\medskip

This matter having been settled, the remainder of the section is
aimed at proving Theorem \ref{thm:mainE}. The following argument
will show that it suffices to prove Theorem \ref{thm:mainE} for {\em
uniformly bounded} processes $\phi$. To see why, for $n\ge 1$ define
\[\phi_n := \one_{\{\n \phi\n \le n\}}\phi.
\]
The processes $\phi_n$ are uniformly bounded, strongly measurable
and adapted, and we have $\limn \phi_n = \phi$ pointwise, hence also
in measure.

We claim that each $\phi_n$ is stochastically integrable with
respect to $W$ and
\[ \limn \int_0^\cdot \phi_n\,dW =\int_0^\cdot \phi\,dW \ \ \hbox{in
measure in $C([0,T];E)$.}
\]
To see this, let $X:\O\to \gg$ be the element represented by $\phi$.
Put
\[
X_n(\omega)f := X(\omega)(\one_{\{\n \phi(\cdot,\omega)\n \le
n\}}f), \qquad f\in\L.
\]
Then by \cite[Proposition 2.4]{NVW},
%$X_n(\omega)\in \gg$, $$\n X_n(\omega)\n_{\gg} \le \n X(\omega)\n_{\gg},$$ and
$\limn X_n = X$ almost surely in $\gg$. It is easily checked that
$\phi_n$ represents $X_n$, and therefore $\phi_n$ is stochastically
integrable by Proposition \ref{prop:NVW}. The convergence in measure
of the stochastic integrals now follows from the continuity
assertion in Proposition \ref{prop:NVW}. This completes the  proof
of  the claim. A more general result in this spirit will be proved
in Section \ref{sec:domination}.

\medskip
It remains to prove Theorem \ref{thm:mainE} for uniformly bounded
processes $\Phi$.

Let $\D_n$ denote the finite $\sigma$-field generated by the $n$-th
dyadic equipartition of the interval $[0,T]$ and let $\G_n =
\D_n\otimes \F$ be the product $\sigma$-field in $[0,T]\times\O$.
Then $\mathbb{G} = \{\G_n\}_{n\ge 1}$ is a filtration in
$[0,T]\times\O$ with $\bigvee_{n\ge 1} \G_n = \mathcal{B}\otimes\F$,
where $\mathcal{B}$ is the Borel $\sigma$-algebra of $[0,T]$. In
what follows with think of $[0,T]\times\O$ as probability space with
respect to the product measure $\frac{dt}{T}\otimes \P$. Note that
for all $f\in L^2([0,T]\times\O;E)$, for almost all $\omega\in \O$
we have
$$%\beq\label{eq:cond}
 \E (f|\G_n)(\cdot,\omega) = \E(f(\cdot,\omega)|\D_n)
 \ \ \hbox{in $L^2(0,T;E)$}.
$$%\eeq
Define the operators $G_n$ on $L^2([0,T]\times \O;E)$ by
\[ G_n f := \tau_n \E(f|\G_n),
\]
where $\tau_n$ denotes the right translation operator over $2^{-n}T$
in $L^2([0,T]\times\O;E)$, i.e., $ \tau_n f(t,\omega) =
\one_{[2^{-n}T,T]} f(t-2^{-n}T,\omega).$

\begin{lemma}\label{lem}
Let $\phi:[0,T]\times\O\to E$ be strongly measurable, adapted,
uniformly bounded, and stochastically integrable with respect to
$W$. Then the processes $\phi_n:[0,T]\times\O\to E$ defined by
$\phi_n := G_n\phi$ are strongly measurable, adapted, uniformly
bounded, and stochastically integrable with respect to $W$.
Moreover, $ \limn \phi_n = \phi$ in measure and \beq\label{eq:conv}
\limn \int_0^\cdot \phi_n\,dW = \int_0^\cdot \phi\,dW \ \ \hbox{in
measure in $C([0,T];E)$}. \eeq
\end{lemma}
\begin{proof} First note that each process $\phi_n$ is strongly measurable,
uniformly bounded, strongly measurable and adapted. By the
vector-valued martingale convergence theorem and the strong
continuity of translations we have
\[
\bal \ & \limn\n \phi - \phi_n\n_{L^2([0,T]\times\O;E)}
\\ &\qquad  \le \limn \n \phi - \tau_n\phi\n_{L^2([0,T]\times\O;E)}
+ \limn \n \tau_n\phi - \tau_n\E(\phi|\G_n)\n_{L^2([0,T]\times\O;E)}
\\ &\qquad \le\limn \n \phi - \tau_n\phi\n_{L^2([0,T]\times\O;E)}
+ \limn \n \phi - \E(\phi|\G_n)\n_{L^2([0,T]\times\O;E)} = 0. \eal
\]
It follows that $\limn \phi_n = \phi$ in $L^2([0,T]\times\O;E)$, and
therefore also in measure.

%Let $n\ge 1$ be fixed for the moment. By \eqref{eq:cond}, for almost
%all $\o\in\O$  the following inequalities hold for all $x\s\in E\s$:
%\beq\label{dom-weak}
%\bal  \int_0^T\!\lb\phi_n(t,\o), x\s\rb^2\,dt
%& \le   \int_0^T\!\lb \E(\phi|\G_n)(t,\o), x\s\rb^2\,dt
%\\ & =  \int_0^T\!\lb \E (\phi(\cdot,\o)|\D_n), x\s\rb^2(t)\,dt
%  =  \int_0^T\!(\E\lb \phi(\cdot,\o), x\s\rb | \D_n)^2(t)\,dt
%\\ & \le  \int_0^T\!\E(\lb \phi(\cdot,\o), x\s\rb^2 | \D_n)\,dt
%  =  \int_0^T\!\lb \phi(t,\o), x\s\rb^2 \,dt.
%\eal\eeq
Let $X:\O\to\gg$ be the random variable represented by $\phi$. For
all $n\ge 1$ let the random variable $X_n:\O\to \gg$ defined by
\[X_n(\omega) := X(\omega)\circ \tau^*_n\circ \E(\,\cdot\,|\D_n)\]
where $\tau_n^*\in \calL(L^2(0,T))$ denotes the left translation
operator. It is easily seen that for all $n\geq 1$, $X_n$ is
represented by $\phi_n$,
%By Lemma \ref{lem:ae-repr},
%for almost all $\o\in \O$ the trajectory $\phi(\cdot,\o)$ represents the
%element $X(\o)\in \gg$.
%Hence by Proposition \ref{prop:domination}, for almost all $\o\in \O$ the
%trajectory $\phi_n(\cdot,\o)$ represents
%an element $X_n(\o)\in \gg$
%. satisfying
%\beq\label{dom-ae}
% \n X_n(\o)\n_{\gg} \le  \n X(\o)\n_{\gg}.
%\eeq
%It is easy to see that
%$\o\mapsto X_n(\o)$ is
%which is easily seen to be strongly measurable, cf. \cite[Lemma 2.5]{NVW}.
%, and it follows from \eqref{dom-ae} that $X_n\in L^p(\O;\gg)$.
%By another application of Lemma \ref{lem:ae-repr},
%$X_n$ is represented by $\phi_n$.
and therefore $\phi_n$ is stochastically integrable with respect to
$W$ by Proposition \ref{prop:NVW}.
%For the proof of \eqref{eq:conv}  we
%fix an arbitrary subsequence $(\phi_{n_k})_{k\ge 1}$ of $(\phi_{n})_{n\ge 1}$.
%Identifying $L^2([0,T]\times \O;E)$ with
%$L^2(\O;L^2(0,T;E))$ we may pass to a further subsequence
%which satisfies $\limj \phi_{{n_k}_j}=  \phi$ in $L^2(0,T;E)$
%almost surely.
By \cite[Proposition 2.4]{NVW}
%and the dominated convergence theorem,
we obtain $ \limj X_{n} = X$ almost surely in $\gg$. Hence, $\limn
X_n = X$ in measure in $\gg$, and \eqref{eq:conv} follows from the
continuity assertion in Proposition \ref{prop:NVW}.
\end{proof}

Now we can complete the proof of Theorem \ref{thm:mainE} for
uniformly bounded processes $\phi$. The processes $\phi_n$ in Lemma
\ref{lem} can be represented as
\[ \phi_n = \sum_{j=1}^{2^n} \one_{I_j} \phi_{j,n},
\]
where the $I_j$ is the $j$-th interval in the $n$-th dyadic
partition of $[0,T]$ and the random variable $\phi_{j,n}:\O\to E$ is
uniformly bounded and $\F_j$-measurable, where $\F_j =\F_{
2^{-n}(j-1)T}$. The proof is completed by approximating the
$\phi_{j,n}$ in $L^0(\O,\F_j;E)$ with simple random variables.
\medskip

Let $E_1$ and $E_2$ be real Banach spaces. Theorem \ref{thm:main}
can be strengthened for $\calL(E_1, E_2)$-valued processes which are
integrable with respect to an $E_1$-valued Brownian motion.

Let $\mu$ be a centred Gaussian Radon measure on $E_1$ and let
$W_\mu$ be an $E_1$-valued Brownian motion with distribution $\mu$,
i.e., for all $t\ge 0$ and $x\s\in E_1\s$ we have
\[
\E \lb W_\mu(t),x\s\rb^2 = t \int_{E_1} \lb x,x\s\rb^2\,d\mu(x).
\]
Let $H_\mu$ denote the reproducing kernel Hilbert space associated
with $\mu$ and let $i_\mu: H_\mu\embed E_1$ be the inclusion
operator. We can associate an $H_\mu$-cylindrical Brownian motion
$W_{H_\mu}$ with $W_\mu$ by the formula
\[
W_{H_\mu}(t) i_\mu\s x\s := \lb W_\mu(t),x\s\rb.
\]

We say $\Phi:[0,T]\times \O\to \calL(E_1,E_2)$ is {\em
$E_1$-strongly measurable and adapted} if for all $x\in E_1$, $\Phi
x$ is strongly measurable and adapted. An $E_1$-strongly measurable
and adapted process $\Phi:[0,T]\times\O\to\calL(E_1,E_2)$ is called
{\em stochastically integrable} with respect to the $E_1$-valued
Brownian motion $W_\mu$ if the process $\Phi\circ i_\mu
:[0,T]\times\O\to\calL(H_\mu, E_2)$ is stochastically integrable
with respect to $W_{H_\mu}$. In this case we write
\[ \int_0^\cdot
\Phi\,dW_\mu := \int_0^\cdot \Phi\circ i_\mu \,dW_{H_\mu}.
\]
By the Pettis measurability theorem and the separability of $H_\mu$,
the $E_1$-strong measurability of $\Phi$ implies the $H_\mu$-strong
measurability of $\Phi\circ i_\mu$. We call $\Phi$ an {\em
elementary adapted} process if $\Phi\circ i_\mu$ is elementary
adapted.

\begin{theorem}\label{thm:main2} Let $E_1$ be a Banach space and let $E_2$ be a UMD space
and fix $p\in (1,\infty)$. Let $W_\mu$ be an $E_1$-valued Brownian
motion with distribution $\mu$. If the process
$\Phi:[0,T]\times\O\to \calL(E_1, E_2)$ is $E_1$-strongly measurable
and adapted and stochastically integrable with respect to $W_\mu$,
there exists a sequence of elementary adapted processes
$\Phi_n:[0,T]\times\O\to \calL(E_1, E_2)$ such that \ben
\item[(i)$''$]
$\displaystyle\limn  \Phi_n x= \Phi x$ in measure for $\mu$-almost
all $x\in E_1$;
\item[(ii)]
$\displaystyle\limn \int_0^\cdot \Phi_n\,dW_\mu =\int_0^\cdot
\Phi\,dW_\mu $ in measure in $C([0,T];E_2)$. \een
\end{theorem}

The proof depends on some well known facts about measurable linear
extensions. We refer to \cite{Bo,FP} for more details. If $\mu$ is a
centred Gaussian Radon measure on $E_1$ with reproducing kernel
Hilbert space $H_\mu$ and $(h_n)_{n\ge 1}$ is an orthonormal basis
$(h_n)_{n\ge 1}$ for $H_\mu$, then the coordinate functionals
$h\mapsto [h,h_n]_{H_\mu}$ can be extended to $\mu$-measurable
linear mappings $\overline{h_n}$ from $E_1$ to $\R$. Moreover, these
extensions are $\mu$-essentially unique in the sense that every two
such extensions agree $\mu$-almost everywhere. Putting
\[\overline{P_n}x := \sum_{j=1}^n
\overline{h_j}x\, h_j, \qquad x\in E_1,
\]
we obtain a $\mu$-measurable linear extension of the orthogonal
projection $P_{n}$ in $H_\mu$ onto the span of the vectors
$h_1,\dots,h_n$. Again this extension is $\mu$-essentially unique,
and we have \beq\label{eq:conv-x}  \limn i_\mu \overline{P_n}x =
\sum_{n\ge 1} \overline{h_j}x\, i_\mu h_j  = x\ \ \hbox{for
$\mu$-almost all $x\in E_1$}. \eeq

\begin{proof}[Proof of Theorem \ref{thm:main2}]
We will reduce the theorem to  Theorem \ref{thm:mainE}. Choose an
orthonormal basis $(h_n)_{n\ge 1}$ of the reproducing kernel Hilbert
space $H_\mu$ and define the processes
$\Psi_n:[0,T]\times\O\to\calL(E_1, E_2)$ by
\[ \Psi_n x :=  \Phi
i_\mu \overline{P_n} x, \qquad x\in E_1.
\]
By \eqref{eq:conv-x},
\[
\limn  \Psi_n x= \Phi x \ \ \hbox{in measure for $\mu$-almost all
$x\in E_1$.}
\]
Also,
\[
\begin{aligned}
\limn \int_0^\cdot \Psi_n\,dW_\mu & =\limn \int_0^\cdot\Psi_{n}\circ
i_\mu\,dW_{H_\mu}
\\ & \stackrel{(*)}{=} \int_0^\cdot\Phi\circ i_\mu\,dW_{H_\mu} = \int_0^\cdot
\Phi\,dW_\mu\ \ \hbox{ in measure in $C([0,T];E_2)$,}
\end{aligned}
\]
where the identity $(*)$ follows by series representation as in the
argument following the statement of Theorem \ref{thm:mainE}. The
proof may now be completed along the lines of this argument; for
$\Phi_k$ we take
\[ \Phi_k x:=
\sum_{n=1}^{N_k} \overline{h_{N_k}}x\, \phi_{j_{k,n},n}, \qquad x\in
E_1,\] where the elementary adapted processes $\phi_{j,n}$
approximate $\Phi i_\mu h_n$ and the indices $N_k$ are chosen as
before.
\begin{comment}
%%
%%
%% the details:
%%
%%
With Theorem \ref{thm:mainE}, for each $n\ge 1$ we choose a sequence
of elementary adapted processes $\phi_{j,n}: [0,T]\times\O\to E_2$
such that $\limj \phi_{j,n} = \Phi i_\mu h_n$ in measure and \[\limj
\int_0^T \phi_{j,n}\,dW_{H_\mu} h_n =\int_0^T \Phi i_\mu
h_n\,dW_{H_\mu} h_n \ \ \hbox{in $L^p(\O;E_2)$.}
\]
Given $k\ge 1$, choose $N_k\ge 1$ so large that
\[ \Big(\E \Big\n \int_0^T
\Phi-\Psi_{N_k}\,dW_\mu\Big\n^p\Bigr)^\frac1p < \frac1{k}.
\]
For each $n=1,\dots,N_k$ choose $j_{k,n}\ge 1$ so large that
\[ \Big|\Big\{\n
\Phi i_\mu h_n - \phi_{j_{k,n},n}\n>\frac1{kN_k}\Big\} \Big| <
\frac{1}{kN_k}
\]
and
\[ \Big(\E \Big\n \int_0^T \Phi i_\mu h_n - \phi_{j_{k,n},n}\,dW_{H_\mu}
h_n\Big\n^p\Bigr)^\frac1p < \frac1{kN_k}.
\]
Define $\Phi_k:[0,T]\times\O\to \calL(E_1, E_2)$ by
\[ \Phi_k x:=  \sum_{n=1}^{N_k}
\overline{h_{N_k}}x\, \phi_{j_{k,n},n}, \qquad x\in E_1.
\]
Each $\Phi_k$ is elementary adapted. As in the proof of Theorem
\ref{thm:main} and using \eqref{eq:conv-x}, it may be shown that
 $\limk \Phi_k x = \Phi x$ in measure for $\mu$-almost all $x\in E_1$
and $\displaystyle \limk\int_0^T \Phi_k\,dW_H=\int_0^T \Phi\,dW_H$
in $L^p(\O;E_2)$. Thus the processes $\Phi_k$ have the desired
properties. %%
%%
%%
\end{comment}
\end{proof}

As a final comment we note that $L^p$-versions of the results of
this section hold as well; for these one has to replace almost sure
convergence by $L^p$-convergence in the proofs.

\section{Domination}\label{sec:domination}

In this section we present two domination results which were
implicit in the arguments so far, and indeed some simple special
cases of them have already been used.

\begin{comment}
The following result follows from a standard application of
Anderson's inequality; cf. \cite{NW} for a discussion.

{\tt $H$-waardig!!}

\begin{proposition}\label{prop:domination}
Let $\phi:[0,T]\to E$ be a function which represents an element
$R\in \gg$. Let $\Psi$ be a set consisting of functions
$\psi:[0,T]\to E$  which satisfy $\lb \psi,x\s\rb\in\L$ and
\[ \int_0^T \lb \psi(t),x\s\rb^2\,dt\le \int_0^T \lb
\phi(t),x\s\rb^2\,dt \ \ \hbox{for all $x\s\in E\s$.}
\]
Then each $\psi\in\Psi$ represents an element $R_\psi \in \gg$ and
we have
\[ \n R_\psi\n_{\gg} \le  \n R\n_{\gg}.
\]
Moreover, the family $\{R_\psi: \ \psi\in\Psi\}$ is relatively
compact in $\gg$.
\end{proposition}

We will use this proposition in the following form. Suppose
$\phi:[0,T]\to E$ represents an element $R\in \gg$ and let
$\phi_n:[0,T]\to E$ be functions which satisfy $\lb
\phi_n,x\s\rb\in\L$ and
\[ \int_0^T \lb \phi_n(t),x\s\rb^2\,dt\le
\int_0^T \lb \phi(t),x\s\rb^2\,dt \ \ \hbox{for all $x\s\in E\s$.}
\]
If $\limn \lb \phi_n, x\s\rb =\lb \phi,x\s\rb$ in $\L$ for all
$x\s\in E\s$, then $\limn R_{\phi_n}=R$ in $\gg$.
\end{comment}

The first comparison result extends \cite[Corollary 4.4]{NW}, where
the case of functions was considered.

\begin{theorem}[Domination]\label{thm:dom} Let
$E$ be a UMD space. Let $\Phi, \Psi :[0,T]\times\O\to\calL(H,E)$ be
$H$-strongly measurable and adapted processes and assume that $\Psi$
is stochastically integrable with respect to $W_H$. If for all
$x\s\in E\s$ we have
$$%\begin{equation}\label{eq:dom}
\int_0^T \n\Phi\s(t)x\s\n_H^2\,dt \le \int_0^T
\n\Psi\s(t)x\s\n_H^2\,dt \ \ \hbox{almost surely},
$$%\end{equation}
then $\Phi$ is stochastically integrable and for all $p\in (1,
\infty)$,
\[\E\sup_{t\in [0,T]}\Bigl\n\int_0^t \Phi(s)\,dW_H(s)\Bigr\n^p \lesssim_{p,E} \E\sup_{t\in [0,T]}\Bigl\n\int_0^t \Psi(s)\,dW_H(s)\Bigr\n^p,\]
whenever the right hand side is finite.
\end{theorem}
\begin{proof}
Since $\Phi$ and $\Psi$ are $H$-strongly measurable and adapted,
without loss of generality we may assume that $E$ is separable.

By Proposition \ref{prop:NVW}, $\Psi$ represents a random variable
$Y:\O\to \ggH$. In particular, for all $x\s\in E\s$ we have $\Psi\s
x\s\in L^2(0,T;H)$ almost surely. We claim that almost surely,
\[
\int_0^T\n\Phi\s(t)x\s\n_H^2\,dt \le  \int_0^T
\n\Psi\s(t)x\s\n_H^2\,dt \ \ \hbox{for all $x\s\in E\s$}.
\]
Indeed, by the reflexivity and separability of $E$ we may choose a
countable, norm dense, $\Q$-linear subspace $F$ of $E^*$. Let $N_1$
be a null set such that \beq\label{eq:domH} \int_0^T
\n\Phi\s(t,\omega)x\s\n_H^2\,dt \le  \int_0^T
\n\Psi\s(t,\omega)x\s\n_H^2\,dt \eeq for all $\omega\in\complement
N_1$ and all $x\s\in F$. By Lemma \ref{lem:ae-repr} there exists a
null set $N_2$ such that $\Psi(\cdot,\omega)$ represents $Y(\omega)$
for all $\omega\in\complement N_2$. Fix $y^*\in E^*$ arbitrary and
choose a sequence $(y_n^*)_{n\geq 1}$ in $F$ such that $\limn y^*_n
= y^*$ in $E\s$ strongly. Fix an arbitrary $\omega\in \complement
(N_1\cup N_2)$. We will prove the claim by showing that
\beq\label{eq:dom2H}\int_0^T \n\Phi\s(t,\omega)y\s\n_H^2\,dt \le
\int_0^T \n\Psi\s(t,\omega)y\s\n_H^2\,dt, \eeq By the closed graph
theorem there exists a constant $C_\omega$ such that
\[\|\Psi^*(\cdot,\omega)x^*\|_{L^2(0,T;H)}\leq C_\omega\|x^*\| \ \ \hbox{for all}  \ x\s\in
E\s.
\]
Hence, $\Psi^*(\cdot,\omega)y^* = \limn \Psi^*(\cdot,\omega)y_n^*$
in $\LH$, by the strong convergence of the $y_n^*$'s to $y\s$. It
follows from \eqref{eq:domH}, applied to the functionals $y_n^*-
y_m^* \in F$, that $(\Phi^* y_n^*)_{n\geq 1}$ is a Cauchy sequence
in $L^2(0,T;H)$. Identification of the limit shows that
$\Phi^*(\cdot,\omega)y^* = \limn \Phi^*(\cdot,\omega)y_n^*$ in
$L^2(0,T;H)$. Now \eqref{eq:dom2H} follows from the corresponding
inequality for $y_n^*$ by letting $n\to \infty$.

By the claim and \cite[Theorem 4.2 and Corollary 4.4]{NW}, almost
every function $\Phi(\cdot,\omega)$ represents an element
$X(\omega)\in \g(\LH,E)$ for which we have
\[\n X(\omega)\n_{\g(\LH,E)}\le \n Y(\omega)\n_{\g(\LH,E)}.\]
By \cite[Remark 2.8]{NVW}) $X$ is strongly measurable as a
$\ggH$-valued random variable. Since $\Phi$ represents $X$, $\Phi$
is stochastically integrable by Proposition \ref{prop:NVW}.
Moreover, from Proposition \ref{prop:NVW} we deduce that
\[\begin{aligned}
\E\sup_{t\in [0,T]}\Bigl\n\int_0^t \Phi(s)\,dW_H(s)\Bigr\n^p
&\eqsim_{p,E} \E\n X\n_{\ggH}^p\leq \E\n Y\n_{\ggH}^p \\ &
\eqsim_{p,E} \E\sup_{t\in [0,T]}\Bigl\n\int_0^t
\Psi(s)\,dW_H(s)\Bigr\n^p.
\end{aligned}\]
\end{proof}

The next result extends \cite[Theorem 6.2]{NW}.
\begin{corollary}[Dominated convergence]\label{cor:DCT}
Let $E$ be a UMD space and fix $p\in (1,\infty)$. For $n\geq 1$, let
$\Phi_n:[0,T]\times\O\to\calL(H,E)$ be $H$-strongly measurable and
adapted and stochastically integrable processes and assume that
there exists an $H$-strongly measurable and adapted process
$\Phi:[0,T]\times\O\to\calL(H,E)$ such that for all $x\s\in E\s$,
\begin{equation}\label{eq:L2conv}
\lim_{n\to\infty} \Phi_n\s x\s = \Phi\s x\s \ \ \hbox{almost surely
in $\LH$.}
\end{equation}
Assume further that there exists an $H$-strongly measurable and
adapted process $\Psi:[0,T]\times\O\to\calL(H,E)$ that is
stochastically integrable and for all $n$ and all $x\s\in E\s$,
\begin{equation}\label{eq:domin}
\int_0^T\n \Phi_n\s (t)x\s\n_{H}^2\,dt \le \int_0^T\n \Psi\s(t)
x\s\n_{H}^2\,dt
 \ \ \hbox{almost surely}.
\end{equation}
Then $\Phi$ is stochastically integrable and
$$
 \limn \int^\cdot_0 \Phi_n- \Phi\,dW_H = 0 \ \ \hbox{in measure in $C([0,T];E)$}.
$$%\end{equation}
\end{corollary}
\begin{proof}
The assumptions \eqref{eq:L2conv} and \eqref{eq:domin} imply that
for all $n$ and $x\s\in E\s$,
\begin{equation}\label{claim}
\int_0^T\n\Phi_n\s (t)x\s\n_{H}^2\,dt \le
\int_0^T\n\Psi\s(t)x\s\n_{H}^2\,dt
 \ \ \hbox{almost surely}.
\end{equation}
Theorem \ref{thm:dom} therefore implies that each $\Phi_n$ is
stochastically integrable, and by passing to the limit $n\to\infty$
in \eqref{claim} we see that the same is true for $\Phi$. Let $Z_n:
\O\to \ggH$ be the element represented by $\Phi_n-\Phi$. By
Proposition \ref{prop:NVW} it suffices to prove that \beq\label{DC3}
\limn Z_n = 0 \ \ \hbox{in measure in $\ggH$.} \eeq As in the proof
of Theorem \ref{thm:dom}, \eqref{eq:domin} implies that for almost
all $\omega\in \O$, \beq\label{eq:domH1} \int_0^T\n\Phi_n\s
(t,\omega)x\s\n_{H}^2\,dt \le
\int_0^T\n\Psi\s(t,\omega)x\s\n_{H}^2\,dt
 \ \ \hbox{for all $n\ge 1$ and $x\s\in E\s$},
 \eeq
and \beq\label{eq:domH2} \int_0^T\n\Phi\s (t,\omega)x\s\n_{H}^2\,dt
\le \int_0^T\n\Psi\s(t,\omega)x\s\n_{H}^2\,dt
 \ \ \hbox{for all $n\ge 1$ and $x\s\in E\s$}.
 \eeq
Denoting by $Y:\O\to \ggH$ the element represented by $\Psi$, we
obtain that,
%\beq\label{eq:domH2} \int_0^T\n\Phi\s (t,\omega)x\s\n_{H}^2\,dt \le
%\int_0^T\n\Psi\s(t,\omega)x\s\n_{H}^2\,dt
% \ \ \hbox{for all $n\ge 1$ and $x\s\in E\s$}.
% \eeq
for almost all $\omega\in \O$, for all $x^*\in E^*$,
\begin{equation}\label{eq:DOM}
\|Z_n^*(\omega) x^*\|_{L^2(0,T;H)} \leq 2\|Y^*(\omega)x^*\|_{L^2(0,T;H)}    %\n Z_n(\omega)\n_{\ggH} \le \, 2 \n Y(\omega)\n_{\ggH}.
\end{equation}
%We claim that for almost all $\omega\in\O$,
%\beq\label{eq:convH}\lim_{n\to\infty} \Phi_n\s(\cdot,\omega) x\s =
%\Phi\s(\cdot,\omega) x\s \ \ \hbox{in $\LH$ for all $x\s\in E\s$}.
%\eeq
Let $N_1$ be a null set such that \eqref{eq:domH1} and
\eqref{eq:domH2} hold for all $\omega\in\complement N_1$. Then for
all $\omega\in \complement N_1$ there is a constant $C(\omega)$
 such that for all $x^*\in E^*$ and all $n\geq 1$,
\begin{equation}\label{constanteCpsi}
\int_0^T\n\Phi\s (t,\omega)-\Phi_n\s (t,\omega)x\s\n_{H}^2\,dt \le
C^2(\omega) \|x^*\|^2.
\end{equation}
Let $(x_j\s)_{j\ge 1}$ be a dense sequence in $E\s$. By
\eqref{eq:L2conv} we can find a null set $N_2$ such that for all
$\omega\in\complement N_2$ and all $j\geq 1$ we have
\beq\label{eq:convj} \lim_{n\to\infty} \Phi_n\s(\cdot, \omega) x\s_j
= \Phi\s(\cdot, \omega) x\s_j \ \text{ in}\ \LH. \eeq Clearly,
\eqref{constanteCpsi} and \eqref{eq:convj} imply that for all
$\omega\in \complement(N_1\cup N_2)$  we have
\[ \lim_{n\to\infty} \Phi_n\s(\cdot, \omega)
x\s = \Phi\s(\cdot, \omega) x\s \ \text{in} \ \LH \ \text{ for all
$x\s\in E\s$},\] %%
%% is wel duidelijk:
%%
%Choose $\omega\in \O_0\cap\O_1$, $x^*\in E^*$ arbitrary. Let $\e>0$ be
%arbitrary. Choose $x_j^*\in E^*$ such that $\|x^*-x_j^*\|<\e$. Then by
%\eqref{constanteCpsi} and the triangle inequality we obtain
%\[\|\Phi^*(\cdot, \omega)x^* - \Phi^*_n(\cdot, \omega)x^*\|_{\LH} \leq
%C(\omega) \e + \|\Phi^*(\cdot, \omega)x^*_j - \Phi^*_n(\cdot,
%\omega)x^*_j\|_{\LH}.\]
%Since $\limn \|\Phi^*(\cdot, \omega)x^*_j - \Phi^*_n(\cdot,
%\omega)x^*_j\|_{\LH} = 0$ and $\e>0$ was arbitrary, the result follows.
%%
%%
hence for almost all $\omega\in \O$, for all $x^*\in E^*$,
\begin{equation}\label{eq:L2conv2}
\limn Z_n^*(\omega) x^* = 0 \ \text{in $L^2(0,T;H)$}.
\end{equation}
By \eqref{eq:DOM} and \eqref{eq:L2conv2} and a standard tightness
argument as in \cite[Theorem 6.2]{NW} we obtain that for almost all
$\omega\in\O$, $\limn Z_n(\omega) = 0$ in $\ggH$. This gives
\eqref{DC3}.
%Now \eqref{DC3} follows from \eqref{eq:DOM}, \eqref{eq:CONV},
%and the dominated convergence theorem.
\end{proof}

Again we leave it to the reader to formulate the $L^p$-version of
these results.

\section{Smoothness - I}\label{sec:besov}

Extending a result of Rosi\'nski and Suchanecki (who considered the
case $H=\R$), it was shown in \cite{NW} (for arbitrary Banach spaces
$E$ and functions $\Phi$) and \cite{NVW} (for UMD Banach spaces and
processes $\Phi$) that if $E$ is a Banach space with type $2$, then
every $H$-strongly measurable and adapted process
$\Phi:[0,T]\times\O\to \calL(H,E)$ with trajectories in
$L^2(0,T;\g(H,E))$ is stochastically integrable with respect to an
$H$-cylindrical Brownian motion $W_H$. Moreover, for $H=\R$ this
property characterises the spaces $E$ with type $2$. Below (Theorem
\ref{thm:besov}) we shall obtain an extension of this result for
processes in UMD spaces with type $p\in [1,2)$.

The results will be formulated in terms of vector valued Besov
spaces. We briefly recall the definition. We follow the approach of
Peetre; see \cite[Section 2.3.2]{Tr} (where the scalar-valued case
is considered) and \cite{Am,GW1,Schm}. The Fourier transform of a
function $f\in L^1(\R^d;E)$ will be normalized  as
$$ \wh{f}(\xi) = \frac1{(2\pi)^{d/2}}\int_{\R^d} f(x)e^{-ix\cdot\xi}\,dx, \quad \xi\in\R^d.$$

Let $\phi\in{\mathscr S}(\R^d)$ be a fixed Schwartz function
%kept fixed in the discussion which follows,
whose Fourier transform $\wh\phi$ is nonnegative and has support in
$\{\xi\in\R^d: \ \tfrac12\le |\xi|\le 2\}$ and which satisfies
$$ \sum_{k\in\Z} \wh\phi(2^{-k}\xi) =1 \quad\hbox{for $\xi\in \R^d\setminus\{0\}$}.$$
Define the sequence $({\varphi_k})_{k\ge 0}$ in ${\mathscr S}(\R^d)$
by
$$\wh{\varphi_k}(\xi) = \wh\phi(2^{-k}\xi) \quad \text{for}\ \  k=1,2,\dots \quad \text{and} \ \ \wh{\varphi_0}(\xi) = 1- \sum_{k\ge 1} \wh{\varphi_k}(\xi), \quad \xi\in\R^d.$$

For $1\le p,q \le \infty$ and $s\in\R$ the {\em Besov space}
$B_{p,q}^s(\R^d;E)$ is defined as the space of all $E$-valued
tempered distributions $f\in {\mathscr S}'(\R^d;E)$ for which
$$ \n f\n_{B_{p,q}^s (\R^d;E)} := \Big\n \big( 2^{ks}{\varphi}_k * f\big)_{k\ge 0} \Big\n_{l^q(L^p(\R^d;E))} $$
is finite. Endowed with this norm, $B_{p,q}^s(\R^d;E)$ is a Banach
space, and up to an equivalent norm this space is independent of the
choice of the initial function  $\phi$. The sequence $({\varphi}_k *
f)_{k\ge 0}$ is called the {\em Littlewood-Paley decomposition} of
$f$ associated with the function $\phi$.
% and we shall write, somewhat informally, $f = \sum_{k\ge 0}{\varphi}_k * f.$

%The following continuous inclusions hold:
%\[\ B_{p, q_1}^s(\R^d; E) \embed B_{p, q_2}^s(\R^d;E), \ B_{p, q}^{s_1}(\R^d;
%E) \embed B_{p, q}^{s_2}(\R^d;E)\] for all $s,s_1, s_2\in \R$, $p, q, q_1,
%q_2\in [1, \infty]$ with $q_1\leq q_2$, $s_2 \leq s_1$. Also note that
%$$ B_{p, 1}^0(\R^d; E)\embed L^p(\R^d;E)\embed B_{p, \infty}^0(\R^d; E).$$
%If $1\le p,q<\infty$, then $B_{p,q}^s(\R^d;E)$ contains the Schwartz space
%${\mathscr S}(\R^d;E)$ as a dense subspace.

Next we define the Besov space for domains. Let $D$ be a nonempty
bounded open domain in $\R^d$. For $1\le p,q\le \infty$ and $s\in\R$
we define
$$B_{p,q}^s (D;E) = \{f|_D: \ f\in B_{p,q}^s (\R^d;E)\}.$$
This space is a Banach space endowed with the norm
$$ \n g\n_{  B_{p,q}^s (D;E)} = \inf_{f|_D = g} \n f\n_{B_{p,q}^s (\R^d;E)}.$$
See \cite[Section 3.2.2]{Tr2} (where the scalar case is considered)
and \cite{Am2}.

We have the following embedding result, which is a straightforward
extension of \cite[Theorems 1.1 and 3.2]{KNVW} where the case $H=\R$
was considered:

\begin{proposition}\label{prop:besov}
Let $E$ be a Banach space and $H$ be a non-zero separable Hilbert
space. Let $D\subseteq \R^d$ be an open domain and let $p\in [1,2]$.
Then $E$ has type $p$ if and only if we have a continuous embedding
\[B_{p,p}^{\frac{d}{p}-\frac{d}{2}}(D;\g(H,E))\embed \g(L^2(D;H);E).\]
\end{proposition}

If we combine this result with Proposition \ref{prop:NVW} we obtain
the following condition for stochastic integrability of processes.

\begin{theorem}\label{thm:besov}
Let $H$ be a separable Hilbert space and let $E$ be a UMD Banach
space with type $p\in [1, 2]$. If $\Phi:[0,T]\times\O\to \calL(H,E)$
is an $H$-strongly measurable process and adapted process with
trajectories in $B_{p,p}^{\frac1p-\frac12}(0,T;\g(H,E))$ almost
surely, then $\Phi$ is stochastically integrable with respect $W_H$.
Moreover, for all $q\in (1, \infty)$,
\[\E\sup_{t\in [0,T]}\Bigl\n\int_0^t \Phi(s)\,dW_H(s)\Bigr\n^q \eqsim_{p,E} \E\|\Phi\|_{B_{p,p}^{\frac1p-\frac12}(0,T;\g(H,E))}^q.\]
\end{theorem}

A similar result can be given for processes with H\"older continuous
trajectories. In particular, invoking \cite[Corollary 3.4]{KNVW} we
see that Theorem \ref{thm:besov} may be applied to functions in
$C^\a([0,1];\g(H,E))$ and, if $E$ is a UMD space, to processes with
paths almost surely in $C^\a([0,1];\g(H,E))$, where
$\a>\frac1p-\frac12$.
%Conversely, by Theorem \ref{thm:counterpart} if for some $\alpha\in
%(0, \frac12)$ all functions in $C^{\a}([0,1];E)$ are stochastically
%integrable with respect to a real Brownian motion, then $E$ has type
%$p$ for all $p\in (1, 2)$ that satisfy $\frac1p-\frac12>\alpha$.
Since UMD spaces always have non-trivial type, %for such $E$
there exists an $\e>0$ such that every $H$-strongly measurable and
adapted process with paths in $C^{\frac{1}{2}-\e}([0,1];\g(H,E))$ is
stochastically integrable with respect to $W_H$. In the converse
direction, \cite[Theorem 3.5]{KNVW} implies that if $E$ is a Banach
space failing type $p\in (1,2)$, then for any $0<\a<\frac1p-\frac12$
there exist examples of functions in $C^\a([0,1];E)$ which fail to
be stochastically integrable with respect to scalar Brownian
motions.

\section{Smoothness - II\label{sec:besov2}}
In this section we give an alternative proof of Proposition
\ref{prop:besov} in the case $D$ is a finite interval. The argument
uses the definition of the Besov space from \cite{Ko} instead of the
Fourier analytic definition of Peetre.

For $s\in (0,1)$ and $p,q\in [1, \infty]$ we will recall the
definition of the Besov space $\Lambda^{s}_{p, q}(0,T;E)$ from
\cite{Ko}. Since it is not obvious that this space is equal to the
Besov space of Section \ref{sec:besov} we use the notation
$\Lambda^{s}_{p, q}(0,T;E)$ instead of $B^{s}_{p, q}(0,T;E)$.

Let $I=(0,T)$. For $h\in\R$ and a function $\phi:I\to E$ we define
the function $T(h)\phi:I\to E$ as the translate of $\phi$ by $h$,
i.e.
\[(T (h) \phi)(t) :=
  \begin{cases}
    \phi(t+h) & \text{if $t+h\in I$}, \\
    0 & \text{otherwise}.
  \end{cases}
\]
For $h\in \R$ put
\[ I[h] : =  \Big\{r\in I:\ r+  h \in I\Big\}.\]
For a strongly measurable function $\phi\in L^p(I;E)$ and $t>0$ let
\[\varrho_p(\phi,t) := \sup_{|h|\le t} \Big(\int_{I[h]} \|T(h) \phi(r) - \phi(r)\|^p \, dr\Big)^{\frac1p}.\]
We use the obvious modification if $p=\infty$.

Now define
\[\Lambda^s_{p, q}(I;\calL(E,F)) := \{\phi\in L^p(I;E): \|\phi\|_{\Lambda^s_{p, q}(I;E)}<\infty\},\]
where
\begin{equation}\label{besov}
\|\phi\|_{\Lambda^s_{p, q}(I;E)} = \Big(\int_0^T \|\phi(t)\|^p \,
dt\Big)^{\frac1p} + \Big(\int_0^1
\big(t^{-s}\varrho_p(\phi,t)\big)^q\, \frac{dt}{t}\Big)^\frac1q
\end{equation}
with the obvious modification for $q=\infty$. Endowed with the norm
$\|\cdot\|_{\Lambda^s_{p, q}(I;E)}$, $\Lambda^s_{p,q}(I;E)$ is a
Banach space.
%Note that our definition is slightly different
%from the one in \cite{Ko}; there, the integral in \eqref{besov} is taken over
%the unit interval. The advantage of the present definition is that it
%simplifies some estimates later on.
The following continuous inclusions hold:
\[\ \Lambda_{p, q_1}^s(I; E) \embed \Lambda_{p, q_2}^s(I;E), \ \Lambda_{p, q}^{s_1}(I;
E) \embed \Lambda_{p, q}^{s_2}(I;E),\] and
\[\Lambda_{p_1, q}^s(I; E) \embed \Lambda_{p_2, q}^s(I;E)\] for all $s,s_1,
s_2\in (0,1)$, $p, p_1, p_2, q, q_1, q_2\in [1, \infty]$ with $1\leq
p_2\le p_1\leq \infty$, $q_1\leq q_2$, $s_2 \leq s_1$.

For all $p\in[1, \infty)$ we have
%the space $\Lambda_{p, q}^s(I; E)$ is the same space as
%the Besov spaces used in Section \ref{sec:besov}. More precisely for all $p\in
%[1, \infty)$
\[\Lambda^s_{p,q}(I;E) = B^s_{p,q}(I;E)\]
with equivalent norms. Here $B^s_{p,q}(I;E)$ is the space defined in
Section \ref{sec:besov}. Since we could not find a reference for
this, we include the short argument. If $I=\R$ this follows from
\cite[Proposition 3.1]{PeWo} (also see \cite[Theorem 4.3.3]{Schm}).
Therefore, for general $I$ the inclusion "$\supseteq$" follows from
the definitions. For the other inclusion notice that by
\cite[Theorem 3.b.7]{Ko} one has
\[\Lambda^s_{p,q}(I;E) = (L^p(I;E), W^{1,p}(I;E))_{s,q}.\]
It is well-known that there is a common extension operator from the
spaces $L^p(I;E)$ and $W^{1,p}(I;E)$ into $L^p(\R;E)$ and
$W^{1,p}(\R;E)$ for all
$p\in [1, \infty]$. Therefore, by interpolation % by\cite[Theorem 1.2.4]{Tr1}
we obtain an extension operator from $(L^p(I;E),
W^{1,p}(I;E))_{s,q}$ into $(L^p(\R;E), W^{1,p}(\R;E))_{s,q}$. Now
the latter is again equal to $B^{s}_{p,q}(\R;E)$ and therefore
"$\subseteq$" holds as well.

We put, for $t>0$,
$$\varphi_p^s(\phi,t):= t^{-s}\varrho_p(\phi,t)
$$
and observe for later use the easy fact that there is a constant
$c_{q,s}>0$ such that for all $\phi\in \Gamma^s_{p, q}(I;E)$ we have
\begin{equation}\label{eq:estBesovsum}
c_{q,s}^{-1} \n\varphi_p^s(\phi,\cdot)\n_{L^q(0,1;\frac{dt}{t})} \le
\big\n \big(\varphi_p^s(\phi,2^{-n})\big)_{n\ge 0}\big\n_{l^q} \le
c_{q,s} \n \varphi_p^s(\phi,\cdot)\n_{L^q(0,1;\frac{dt}{t})}.
\end{equation}

\begin{theorem}\label{thm:BesovembM}
Let $H$ be a separable Hilbert space, $E$ a Banach space, and let
$p\in [1,2)$. Then $E$ has type $p$ if and only if $\Lambda_{p,
p}^{\frac1p-\frac12}(0,T;\g(H,E)) \hookrightarrow \g(L^2(0,T;H),E)$
continuously.
\end{theorem}

\begin{proof}
For the proof that $E$ has type $p$ if the inclusion holds we refer
to \cite[Theorem 3.3]{KNVW}. To prove the converse we may assume
$T=1$. Let $(g_{00}, g_{nk}:n\geq 0, k = 1, \ldots, 2^n)$ be the
$L^2$-normalized Haar system on $[0,1]$, i.e. $g_{00}\equiv1$ and
for all other $n$ and $k$ let
\[g_{nk} =
  \begin{cases}
    2^{\frac{n}{2}} & \text{on}\ [(k-1)2^{-n}, (k-\frac12)2^{-n}) \\
    -2^{\frac{n}{2}} & \text{on}\ [(k-\frac12)2^{-n}, k2^{-n}) \\
    0 & \text{otherwise}.
  \end{cases}
\]
Let $(h_i)_{i\geq 1}$ be an orthonormal basis for $H$. Then
$(g_{nk}\otimes h_i)_{m,k,i}$ is an orthonormal basis for
$L^2(0,1;H)$. Let $(\g_i)_i$, $(\g_{nki})_{n,k,i}$ be Gaussian
sequences and let $(r_{nk})_{n,k}$ be an independent Rademacher
sequence. Let $\Phi\in \Lambda_{p,
p}^{\frac1p-\frac12}(0,T;\g(H,E))$ be arbitrary. Since $E$ has type
$p$, $L^2(\O;E)$ has type $p$ with $T_p(L^2(\O;E))=T_p(E)$ (cf.
\cite{DJT}) and we have
\[\begin{aligned}
\Big(&\E\Big\|\sum_{i\geq 1} \g_{nki} I_\Phi g_{00} \otimes h_i +
\sum_{n\geq 0} \sum_{k=1}^{2^n} \sum_{i\geq 1} \g_{nki} I_\Phi
g_{nk} \otimes h_i\Big\|^2\Big)^{\frac12} \\ & = \Big(\E_r
\E\Big\|\sum_{i\geq 1} \g_{nki} I_\Phi g_{00} \otimes h_i +
\sum_{n\geq 0} \sum_{k=1}^{2^n} \sum_{i\geq 1} r_{nk} \g_{nki}
I_\Phi g_{nk} \otimes h_i\Big\|^2\Big)^{\frac12}
\\ & \leq \Big\|\sum_{i\geq 1} \g_{i} I_\Phi g_{00} \otimes
h_i\Big\|_{L^2(\O;E)} + T_p(E) \Big(\sum_{n\geq 0} \sum_{k=1}^{2^n}
\Big\|\sum_{i\geq 1} \g_{i} I_\Phi g_{nk} \otimes
h_i\Big\|_{L^2(\O;E)}^p\Big)^{\frac1p}
\end{aligned}\]
Now one easily checks that
\[\Big\|\sum_{i\geq 1} \g_{i} I_\Phi g_{00} \otimes
h_i\Big\|_{L^2(\O;E)}\leq \|\Phi\|_{L^p(0,1;\g(H,E))}.\] For the
other term note that
\[I_\Phi g_{nk} \otimes h_i = 2^{\frac{n}{2}} \int_{(k-1)2^{-n}}^{(k-\frac12)2^{-n}} (\Phi(s) - \Phi(s+2^{-n-1})) h_i \, ds.\]
Therefore,
\[\begin{aligned}
\sum_{k=1}^{2^n}  \Big\|&\sum_{i\geq 1} \g_{i} I_\Phi g_{nk} \otimes
h_i\Big\|^p_{L^2(\O;E)} \\ & = 2^{\frac{np}{2}}  \sum_{k=1}^{2^n}
\Big\| \int_{(k-1)2^{-n}}^{(k-\frac12)2^{-n}} \Phi(s) -
\Phi(s+2^{-n-1}) \, ds \Big\|^p_{\g(H,E)}
\\ & \leq 2^{\frac{np}{2}} 2^{(n+1)(1-p)} \sum_{k=1}^{2^n}
\int_{(k-1)2^{-n}}^{(k-\frac12)2^{-n}} \|\Phi(s) -
\Phi(s+2^{-n-1})\|^p_{\g(H,E)} \, ds
\\ & \leq 2^{-p + 1} 2^{n(1-\frac{p}{2})}
\int_{0}^{1-2^{-n-1}} \|\Phi(s) - \Phi(s+2^{-n-1})\|^p_{\g(H,E)} \,
ds
\end{aligned}\]
We conclude that
\[\begin{aligned}
\Big(\sum_{n\geq 0} &\sum_{k=1}^{2^n}  \E\Big\|\sum_{i\geq 1} \g_{i}
I_\Phi g_{nk} \otimes h_i\Big\|^p\Big)^{\frac1p} \\ & \leq 2^{-1 +
\frac1p} \Big(\sum_{n\geq 0} 2^{n(1-\frac{p}{2})}
\int_{0}^{1-2^{-n-1}} \|\Phi(s) - \Phi(s+2^{-n-1})\|^p_{\g(H,E)} \,
ds\Big)^{\frac1p}
\\ & \lesssim_p \|\Phi\|_{\Lambda_{p,
p}^{\frac1p-\frac12}(0,T;\g(H,E))},
\end{aligned}\]
where the last inequality follows from \eqref{eq:estBesovsum}.
\end{proof}

If $(0,T)$ is replaced with $\R$, one can use the Haar basis on each
interval $(j, j+1)$ to obtain the analogous embedding result for
$\R$. More generally, the proof can be adjusted to the case of
finite or infinite rectangles $D\subseteq \R^d$. Furthermore, using
extension operators one can extend the embedding result to bounded
regular domains.

\medskip

As a consequence of Theorem \ref{thm:BesovembM} we recover a
H\"older space embedding result from \cite{KNVW}. Using
\cite[Theorem 2.3]{NW} this can be reformulated as follows.

\begin{proposition}\label{prop:embtype}
Let $E$ be a Banach space and let $p\in [1, 2)$. If $E$ type $p$,
then for all $\alpha>\frac{1}{p}-\frac12$ it holds that $\phi\in
C^\alpha([0,1];E)$ implies that $\phi$ is stochastically integrable
with respect to $W$. Moreover, there exists a constant $C$ only
depending on the type $p$ constant of $E$ such that
\[\E\Big\|\int_0^1 \phi \, d W\Big\|^2\leq C^2\|\phi\|_{C^{\alpha}([0,1];E)}^2\]
\end{proposition}

In \cite{KNVW} a converse to this result is obtained as well: if all
functions in $C^\alpha([0,1];E)$ are stochastically integrable, then
$E$ has type $p$ for all $p\in [1, 2)$ satisfying
$\alpha<\frac1p-\frac12$. However, the case that
$\alpha=\frac1p-\frac12$ is left open there and will be considered
in the following theorem. For the definition of stable type $p$ we
refer to \cite{LeTa}.

\begin{theorem}\label{thm:stabletype}
Let $E$ be a Banach space, let $\alpha\in (0, \frac12]$ and let
$p\in [1,2)$ be such that $\alpha=\frac1p-\frac12$. If every
function in $C^{\alpha}([0,1];E)$ is stochastically integrable with
respect to $W$, then $E$ has stable type $p$.
\end{theorem}

Since $l^p$ spaces for $p\in [1, 2)$ do not have stable type $p$, it
follows from Theorem \ref{thm:stabletype} that there exists a
$(\frac1p-\frac12)$-H\"older continuous function $\phi:[0,1]\to l^p$
that is not stochastically integrable with respect to $W$. An
explicit example can be obtained from the construction below. This
extends certain examples in \cite{RS,Yo}

\begin{proof}

{\em Step 1:} Fix an integer $N\geq 1$. First we construct an
certain function with values in $l^p_N$.
%Let $h_{00}$, $h_{nk}$ for $n\geq 0, k = 1, \ldots, 2^n$ be the
%$L^2$-normalized Haar system on $[0,1]$,
% i.e. $h_{00}\equiv1$ and for all other $n$ and $k$ let
%\[h_{nk} =
%  \begin{cases}
%    2^{\frac{n}{2}} & \text{on}\ [(k-1)2^{-n}, (k-\frac12)2^{-n}) \\
%    -2^{\frac{n}{2}} & \text{on}\ [(k-\frac12)2^{-n}, k2^{-n}) \\
%    0 & \text{otherwise}.
%  \end{cases}
%\]
Let $\varphi_{00}, \varphi_{nk}$ for $n\geq 0, k = 1, \ldots, 2^n$
be the Schauder functions on $[0,1]$, i.e., $\varphi_{nk}(x) =
\int_0^x g_{nk}(t) \, dt$ where $g_{nk}$ are the $L^2$-normalized
Haar functions. Let $(e_{n})_{n=1}^N$ be the standard basis in
$l^p_N$. Let $\psi:[0,1]\to l^p$ be defined as
\[\psi(t) = \sum_{n=0}^N \sum_{k=1}^{2^n} 2^{\frac{(p-1)n}{p}}\varphi_{nk}(t) e_{2^n+k}.\]
Then $\psi$ is stochastically integrable and
\[\begin{aligned}
\E\Big\| \int_0^1 \psi \, d W \Big\|^p &= \E\Big\| \sum_{n=0}^N
\sum_{k=1}^{2^n} 2^{\frac{(p-1)n}{p}} \int_0^1 \varphi_{nk} \, d W
e_{2^n+k} \Big\|^p
\\ & = \sum_{n=0}^N \sum_{k=1}^{2^n} 2^{(p-1)n} \E\Big|\int_0^1 \varphi_{nk} \, d W\Big|^p
\\ & = m_p^p \sum_{n=0}^N \sum_{k=1}^{2^n} 2^{(p-1)n} \Big(\E\Big|\int_0^1 \varphi_{nk} \, d
W\Big|^2\Big)^{\frac{p}{2}}
\\ & = m_p^p \sum_{n=0}^N \sum_{k=1}^{2^n} 2^{(p-1)n} \Big(\int_0^1 \varphi_{nk}^2(t) \, d
t\Big)^{\frac{p}{2}}
\\ & = m_p^p \sum_{n=0}^N \sum_{k=1}^{2^n} 2^{(p-1)n} \Big( \frac{2^{-2n-2}}{3} \Big)^{\frac{p}{2}} = m_p^p \frac{N
}{12^{\frac{p}{2}}},
\end{aligned}\]
where $m_p = (\E |W(1)|^p)^{\frac1p}$. Therefore,
\begin{equation}\label{eq:stochint}
\Big(\E\Big\| \int_0^1 \psi \, d W \Big\|^2\Big)^{\frac12} \geq
\Big(\E\Big\| \int_0^1 \psi \, d W \Big\|^p\Big)^{\frac1p} = K_p
N^{\frac1p},
\end{equation}
with $K_p = m_p/\sqrt{12}$.

On the other hand $\psi$ is $\alpha$-H\"older continuous with
\begin{equation}\label{eq:holdernorm}
\|\psi\|_{C^{\alpha}([0,1];E)} = \sup_{t\in [0,1]}\|\psi(t)\| +
\sup_{0\leq s<t\leq 1} \frac{\|\psi(t) - \psi(s)\|}{(t-s)^{\alpha}}
\leq C_p,
\end{equation} where $C_p$ is a constant only depending on
$p$. Indeed, for each $t\in [0,1]$, we have
\[\begin{aligned}
\|\psi(t)\| &= \Big(\sum_{n=0}^N \sum_{k=1}^{2^n} 2^{(p-1)n}
|\varphi_{nk}(t)|^p\Big)^{\frac1p} \leq \Big(\sum_{n=0}^N
 2^{(p-1)n} 2^{-\frac{(n+2)p}{2}}\Big)^{\frac1p}
\\ & = \Big(\sum_{n=0}^N 2^{-(1-p/2)n - p }\Big)^{\frac1p}
\leq \frac12 \Big(\frac{2^{1-p/2}}{2^{1-p/2}-1}\Big)^{\frac1p}.
\end{aligned}\]
Now fix $0\leq s<t\leq 1$. Let $n_0$ be the largest integer such
that there exists an integer $k$ with the property that $s,t\in
[(k-1)2^{-n_0}, (k+1) 2^{-n_0}]$. Then\footnote{This argument
corrects a minor mistake in the proof in the published version of
the paper.}
\[(k-1)2^{-n_0} \leq s\leq k2^{-n_0}\leq t\leq (k+1)2^{-n_0}.\]
Indeed, otherwise one can replace $n_0$ with $n_0+1$ using the point
$(k-1/2) 2^{-n_0} = (2k-1) 2^{-(n_0+1)}$ or $(k+1/2) 2^{-n_0} =
(2k+1) 2^{-(n_0+1)}$ as the middle of the dyadic interval with
length $2^{-(n_0+1)}$. Moreover, $2^{-(n_0+1)}\leq (t-s)\leq
2^{-n_0+1}$. The upper estimate is clear. For the lower estimate
assume that $(t-s)<2^{-(n_0+1)}$. Then $(t-k 2^{-n}) < 2^{-(n_0+1)}$
and $(k 2^{-n}-s) < 2^{-(n_0+1)}$. Therefore, $t,s\in
[(k-1/2)2^{-n_0}, (k+1/2) 2^{-n_0}] = [(2k-1)2^{-(n_0+1)}, (2k+1)
2^{-(n_0+1)}]$. This contradicts the choice of $n_0$.

Now for each $0\leq n\leq n_0$ let $k_n$ be the unique integer such
that $s\in [(k_n-1)2^{-n}, k_n 2^{-n})$. Now two cases occur: (i)
$t\in [(k_n-1)2^{-n}, k_n 2^{-n}]$ or (ii) $t\in [k_n2^{-n}, (k_n+1)
2^{-n}]$.

In case (i) it follows that
\[|\varphi_{n k_n}(t) - \varphi_{nk_n}(s)|\leq 2^{\frac{n}{2}} (t-s) \leq 2^{\frac{n}{2}} 2^{(-n_0+1)(1-\alpha)} (t-s)^\alpha.\]
In case (ii) it follows that
\[\begin{aligned}
|\varphi_{n k_n}(t) - \varphi_{nk_n}(s)|& = |\varphi_{nk_n}(k_n
2^{-n}) - \varphi_{nk_n}(s)| \leq 2^{\frac{n}{2}} (k_n 2^{-n} -  s)
\\ & \leq 2^{\frac{n}{2}} (t -  s) \leq  2^{\frac{n}{2}}
2^{(-n_0+1)(1-\alpha)} (t-s)^\alpha
\end{aligned}\]
and in the same way
\[|\varphi_{n k_n+1}(t) - \varphi_{nk_n+1}(s)| = |\varphi_{nk_n+1}(t)- \varphi_{nk_n+1}(k_n 2^{-n})| \leq  2^{\frac{n}{2}} 2^{(-n_0+1)(1-\alpha)} (t-s)^\alpha.\]

For $n_0< n\leq N$ let $k_n>\ell_n$ be the unique integers such that
$t\in [(k_n-1)2^{-n}, k_n 2^{-n}]$ and $s\in [(\ell_n-1)2^{-n},
\ell_n 2^{-n}]$. Then
\[|\varphi_{n k_n}(t) - \varphi_{nk_n}(s)| = |\varphi_{n k_n}(t)|\leq 2^{-\frac{n}{2}-1},\]
\[|\varphi_{n \ell_n}(t) - \varphi_{n \ell_n}(s)| = |\varphi_{n \ell_n}(s)|\leq 2^{-\frac{n}{2}-1}.\]
We conclude that
\[\begin{aligned}
\|&\psi(t) - \psi(s)\|^p \\ & \ \leq \sum_{n=0}^{n_0}
\sum_{k=1}^{2^n} 2^{(p-1)n} |\varphi_{nk}(t) - \varphi_{nk}(s)|^p +
\sum_{n=n_0+1}^{N} \sum_{k=1}^{2^n} 2^{(p-1)n} |\varphi_{nk}(t) -
\varphi_{nk}(s)|^p
\\ & \leq 2\sum_{n=0}^{n_0}  2^{(p-1)n} 2^{\frac{n p }{2}} 2^{(-n_0+1)(1-\alpha)p}
(t-s)^{\alpha p} + \sum_{n=n_0+1}^{N} 2^{(p-1)n} 2^{-\frac{n p}{2}}
\\ & \leq \frac{2^{2-\alpha}}{2^{\frac32p -1}-1} (t-s)^{\alpha p} +
\frac{2^{-(n_0+1)(1-\frac{p}{2})}}{1-2^{-(1-\frac{p}{2})}}
\end{aligned}\]
Noting that $2^{-(n_0+1)}\leq (t-s)$ and $(1-\frac{p}{2})=\alpha p$
it follows that
\[\|\psi(t) - \psi(s)\|\leq \Big(\frac{2^{2-\alpha}}{2^{\frac32p -1}-1}+ \frac{1}{1-2^{-(1-\frac{p}{2})}}\Big)^{\frac1p}
(t-s)^{\alpha}.\] Therefore, \eqref{eq:holdernorm} follows.

{\em Step 2:} Assume that every function in $C^{\alpha}([0,1];E)$ is
stochastically integrable. It follows from the closed graph theorem
that there exists a constant $C$ such that for all $\phi\in
C^{\alpha}([0,1];E)$ we have
\begin{equation}\label{eq:estforE}
\Big(\E\Big\|\int_0^1 \phi \, d W\Big\|^2\Big)^{\frac12}\leq
C\|\phi\|_{C^{\alpha}([0,1];E)}
\end{equation}
Now assume that $E$ does not have stable type $p$. By the
Maurey-Pisier theorem \cite[Theorem 9.6]{LeTa} it follows that $l^p$
is finitely representable in $E$. In particular it follows that for
each integer $N$ there exists an operator $T_N:l^p_N\to E$ such that
$\|x\|\leq \|T_N x\|\leq 2\|x\|$ for all $x\in l^p_N$. Now let
$\phi:[0,1]\to E$ be defined as $\phi(t) = T_N \psi_N(t)$, where
$\psi_N:[0,1]\to l^p_N$ is the function constructed in Step 1. Then
it follows from \eqref{eq:stochint}, \eqref{eq:holdernorm} and
\eqref{eq:estforE} that
\[\begin{aligned}
K_p N^{\frac1p} \leq \Big(\E\Big\|\int_0^1 \psi \, d
W\Big\|^2\Big)^{\frac12} &\leq \Big(\E\Big\|\int_0^1 \phi \, d
W\Big\|^2\Big)^{\frac12}
\\ & \leq C\|\phi\|_{C^{\alpha}([0,1];E)} \leq 2 C\|\psi\|_{C^{\alpha}([0,1];l^p_N)}\leq 2 C C_p.
\end{aligned}\]
This cannot hold for $N$ large and therefore $E$ has stable type
$p$.
\end{proof}

As a corollary we obtain that the set of all $\alpha\in (0,
\frac12]$ such that every $f\in C^{\alpha}([0,1];E)$ is
stochastically integrable is relatively open.
\begin{corollary}
Let $E$ be a Banach space and let $\alpha\in (0, \frac12]$ and let
$p\in [1,2)$ be such that $\alpha=\frac1p-\frac12$. If every
function in $C^{\alpha}([0,1];E)$ is stochastically integrable with
respect to $W$, then $E$ has (stable) type $p_1$ for some $p_1>p$.
In particular, there exists an $\e\in (0,\alpha)$ such that every
function in $C^{\alpha-\e}([0,1];E)$ is stochastically integrable.
\end{corollary}
\begin{proof}
The first part follows from Theorem \ref{thm:stabletype} and
\cite[Corollary 9.7, Proposition 9.12]{LeTa}. The last statement is
a consequence of this and Proposition \ref{prop:embtype}, where
$\e>0$ may be taken such that $\alpha-\e = \frac{1}{p_1}-\frac12$
\end{proof}

\section{ Banach function spaces}\label{sec:BFS}

In this section we prove a criterium (Theorem \ref{thm:funcspace})
for stochastic integrability of a process in the case $E$ is a UMD
Banach function space which was stated without proof in \cite{NVW}.
It applies to the spaces $E=L^p(S)$, where $p\in (1, \infty)$ and
$(S,\Sigma,\mu)$ is a $\sigma$-finite measure space.

 We start with the case where
$\Phi$ is a function with values in $\calL(H,E)$. The following
proposition extends \cite[Corollary 2.10]{NW}, where the case $H=\R$
was considered.

\begin{proposition}\label{fs}
Let $E$ be Banach function space with finite cotype over a
$\sigma$-finite measure space $(S,\Sigma,\mu)$. Let $\Phi:[0,T]\to
\calL(H,E)$ be an $H$-strongly measurable function and assume that
there exists a strongly measurable function $\phi:[0,T]\times S\to
H$ such that for all $h\in H$ and $t\in [0,T]$,
\[(\Phi(t) h)(\cdot) = [\phi(t,\cdot), h]_H \ \ \hbox{in $E$}.\]
Then $\Phi$ is stochastically integrable if and only if
\begin{equation}\label{eq:fs}
 \Big\n \Big(\int^T_0 \|\phi(t,\cdot)\|_H^2\,dt\Big)^\frac12\Big\n_E<\infty.
\end{equation}
In this case we have
\[
\Big( \E  \Big\n \int_0^T \Phi\,d W_H\Big\n_E^2\Big)^\frac12
 \eqsim_E
\Big\n \Big(\int^T_0 \|\phi(t,\cdot)\|_H^2\,dt\Big)^\frac12\Big\n_E.
\]
\end{proposition}
\begin{proof}
First assume that $\Phi$ is stochastically integrable. Let
$\mathcal{N} = \{n\in \mathbb{N}: 1\leq n< {\rm dim}(H) +1\}$, let
$(e_m)_{m\in\mathcal{N}}$ be the standard unit basis for
$L^2(\mathcal{N},\tau)$, where $\tau$ denotes the counting measure
on $\mathcal{N}$. Choose orthonormal bases $(f_n)_{n\ge 1}$ for $\L$
and $(h_n)_{n\in\mathcal{N}}$ for $H$. Define $\Psi:[0,T]\times
\mathcal{N} \to  E$ by $\Psi(t, n) := \Phi(t) h_n$ and define the
integral operator $I_{\Psi}:L^2([0,T]\times
\mathcal{N}\!,\,dt\times\tau) \to E$ by
\[ I_\Psi f := \int_{\mathcal{N}}\int_{[0,T]} f(t,n)\Psi(t,n)\,dt\,d\tau(n)
= \sum_{n\in\mathcal{N}}\int_0^T f(t,n)\Phi(t)h_n\,dt.
\]
Note that the integral on the right-hand side is well defined as a
Pettis integral. Let $I_\Phi\in \g(\LH,E)$ be the operator
representing $\Phi$ as in Proposition \ref{prop:NVW} (the special
case for functions). Then  $I_\Psi\in \g(L^2([0,T]\times
\mathcal{N}\!,\,dt\times\tau),E)$ and
\[\Big( \E  \Big\n \int_0^T \Phi\,d W_H\Big\n_E^2\Big)^\frac12
  = \|I_{\Phi}\|_{\g(\LH, E)}
  = \|I_{\Psi}\|_{\g(L^2([0,T]\times \mathcal{N}\!,\,dt\times\tau),E)}. \]
On the other hand, by a similar calculation as in \cite[Corollary
2.10]{NW} one obtains, with $(r_{mn})$ denoting a doubly indexed
sequence of Rademacher variables on a probability space
$(\O',\F',\P')$,
\[\begin{aligned}
\|I_{\Psi}\|_{\g(L^2([0,T]\times \mathcal{N}\!,\,dt\times\tau),E)}
&\eqsim_E \Big(\E'\Big\|\sum_{m,n} r_{mn} \int_0^T \sum_{k} \Psi(t,
k) e_m(k) f_n(t) \, dt  \Big\|_E^2 \Big)^{\frac12}
\\ &\eqsim_E \Big\n \Big(\int^T_0 \sum_{k}
\big|\Psi(t,k)(\cdot)\big|^2 \,dt \Big)^\frac12\Big\n_E
%\\ & = \Big\n \Big(\int^T_0 \sum_{k}\big|(\Phi(t)
%h_k)(\cdot)\big|^2\,dt\Big)^\frac12\Big\n_E \\ &= \Big\n
%\Big(\int^T_0 \sum_{k}\big|[\phi(t,\cdot),
%h_k]_H\big|^2\,dt\Big)^\frac12\Big\n_E
\\ & = \Big\n \Big(\int^T_0
\|\phi(t,\cdot)\|_H^2\,dt\Big)^\frac12\Big\n_E.
\end{aligned}\]

For the converse one can read all estimates backwards, but we have
to show that $\Phi$ belongs to $\LH$ scalarly if \eqref{eq:fs}
holds. For all $x^*\in E^*$ we have
\[\begin{aligned}
\|\Phi^*x^*\|_{\LH}^2 & = \Big(\sum_{m,n} \Big(\int_0^T [\Phi^*(t)
x^*, h_m]_H f_n(t) \, dt\Big)^2 \Big)^{\frac12}
\\ &
=  \Big(\sum_{n, m} \Big(\int_0^T \sum_{k} \lb \Psi(t, k), x^*\rb
e_m(k) f_n(t) \, dt \Big)^2 \Big)^{\frac12}
%\\ & = \Big(\E'\Big|\sum_{n, m} r_{mn} \int_0^T \sum_{k}  \lb \Psi(t, k), x^*\rb e_m(k) f_n(t) \, dt  \Big|^2 \Big)^{\frac12}
\\ & \leq \Big(\E'\Big\|\sum_{n, m} r_{mn} \int_0^T \sum_{k} \Psi(t, k) e_m(k)
f_n(t) \, dt  \Big\|_E^2 \Big)^{\frac12}\|x^*\|.
\end{aligned}\]
\end{proof}

By combining this proposition with Proposition \ref{prop:NVW} and
recalling the fact that UMD spaces have finite cotype, we obtain:

\begin{theorem}\label{thm:funcspace}
Let $E$ be UMD Banach function space over a $\sigma$-finite measure
space $(S,\Sigma,\mu)$ and let $p\in(1, \infty)$. Let
$\Phi:[0,T]\times\O \to \calL(H,E)$ be an $H$-strongly measurable
and adapted process and assume that there exists a strongly
measurable function $\phi:[0,T]\times \O\times S\to H$ such that for
all $h\in H$ and $t\in [0,T]$,
\[(\Phi(t) h)(\cdot) = [\phi(t,\cdot), h]_H \ \ \hbox{in $E$}.\]
Then $\Phi$ is stochastically integrable if and only if
\[
\Big\n \Big(\int^T_0
\|\phi(t,\cdot)\|_H^2\,dt\Big)^\frac12\Big\n_E<\infty \ \
\hbox{almost surely}.
\]
In this case for all $p\in (1, \infty)$ we have
\[
\E  \sup_{t\in [0,T]}\Big\n \int_0^t \Phi(t)\,d W_H(t)\Big\n^p
\eqsim_{p,E} \E\Big\n \Big(\int^T_0
\|\phi(t,\cdot)\|_H^2\,dt\Big)^\frac12\Big\n_E^p.
\]
\end{theorem}

{\em Acknowledgment} -- The second named author thanks S. Kwapie\'n
for the helpful discussion that led to Theorem \ref{thm:stabletype}.

\providecommand{\bysame}{\leavevmode\hbox
to3em{\hrulefill}\thinspace}

\end{document}